\documentclass[reqno]{amsart}
\usepackage{url}
\newtheorem{theorem}{Theorem}[section]

\newtheorem{conjecture}[theorem]{Conjecture}

\theoremstyle{definition}

\theoremstyle{remark}

\numberwithin{equation}{section}
\usepackage{algorithmic}
\usepackage[top=3cm,bottom=2cm,right=2cm,left=2cm]{geometry}

\begin{document}

\title[The Application of Modular Functions to Infinite Partition Congruence Families]{Divisibility Arising From Addition: The Application of Modular Functions to Infinite Partition Congruence Families}


\author{}
\address{}
\curraddr{}
\email{}
\thanks{}

\author{Nicolas Allen Smoot}
\address{}
\curraddr{}
\email{}
\thanks{}



\date{}

\dedicatory{}

\begin{abstract}
The theory of partition congruences has been a fascinating and difficult subject for over a century now.  In attempting to prove a given congruence family, multiple possible complications include the genus of the underlying modular curve, representation difficulties of the associated sequences of modular functions, and difficulties regarding the piecewise $\ell$-adic convergence of elements of the associated space of modular functions.  However, our knowledge of the subject has developed substantially and continues to develop.  In this very brief survey, we will discuss the utility of modular functions in proving partition congruences, both theoretical and computational, and many of the problems in the subject that are yet to be overcome.
\end{abstract}

\maketitle

\section{Introduction}

A common uninformed criticism of science is that analyzing a given subject detracts from its ``wonder and beauty" (e.g., \cite{Munroe}).  While it is of course true that wonder and beauty are largely subjective concepts, many working scientists might argue that the reverse of this criticism is generally true: studying a subject which at first sight appears banal can often reveal an astonishing hidden structure with interconnections wholly unexpected.

One of the most remarkable examples of this counterargument resides in the study of addition over the whole numbers---what is called the theory of integer partitions.  Few subjects in mathematics are easier to grasp.

We define a partition of a given $n\in\mathbb{N}=\mathbb{Z}_{\ge 1}$ as an expression of $n$ as a sum of other members of $\mathbb{N}$, called parts, in which parts may be repeated, and ordering of the parts is irrelevant.  Thus the number 5 can be partitioned as $5,\ 4+1,\ 3+2,\ 3+1+1,\ 2+2+1,\ 2+1+1+1,\ 1+1+1+1+1$.

The total number of partitions of a given $n$ is denoted $p(n)$, and is also a member of $\mathbb{N}$.  Thus, $p(5)=7$.  For technical reasons, we define $p(0):=1$, and $p(x)=0$ for $x\not\in\mathbb{Z}_{\ge 0}$.

At first sight, this subject appears especially simple.  Ahlgren and Ono have referred to it as ``child's play" \cite{Ahlgren3}, while Hirschhorn has referred to his work in the subject as ``high school algebra, but taken somewhat further" \cite{HirschhornQ}.  To understand it, surely one only needs to understand whole numbers and addition.  We would not expect that, say, functions of a complex variable, or the theory of complex manifolds would play a substantial role in the subject.  On the other hand, the sequence $\left(p(n)\right)_{n\ge 0}$ begins $$(1,1,2,3,5,7,11,15,22,30,42,56,77,101,135,176,231,297,385,490,627,1002,1255,1575,1958,2436,3010,...).$$  Aside from a brief overlap with the Fibonacci numbers, the specific values of the sequence have an apparently pseudorandom appearance.  Indeed, $p(n)$ has been compared to the prime counting function $\pi(x)$ in the past \cite[31:30]{Sykes}, suggesting that the subject would be highly resistant to detailed study.

We briefly consider the question of an \textit{efficient formula} for $p(n)$.  We give the following formula, originally given in an incomplete form by Hardy and Ramanujan in 1918 \cite{HardyR}, and which was refined by Rademacher in 1937 \cite{Rademacher36}.

\begin{align}
p(n) = \frac{1}{\pi\sqrt{2}}\sum_{k=1}^{\infty}\sqrt{k}\sum_{\substack{0\le h<k,\\ \mathrm{gcd}(h,k)=1}}e^{-2\pi i n h/k + \pi i s(h,k)} \frac{d}{dx}\left( \frac{\sinh\left( \frac{\pi}{k}\sqrt{\frac{2}{3}\left(x-\frac{1}{24}\right)} \right)}{\sqrt{x-\frac{1}{24}}} \right)\Bigg\rvert_{x=n},\label{pnformula}
\end{align} in which $s(h,k)\in\mathbb{Q}$ by a certain sawtooth function.  It should be noted that this formula provides a \textit{near-optimal} computational efficiency \cite{Johansson}.

This astonishing formula owes its existence largely to the generating function for $p(n)$, which we define as

\begin{align}
\mathcal{P}(q) := \sum_{n=0}^{\infty}p(n)q^n = \prod_{m=1}^{\infty}\frac{1}{1-q^m}.
\end{align}  This can be shown by performing a Taylor series expansion on each factor $(1-q^m)^{-1}$ and examining the resulting coefficient of $q^n$.  What is remarkable is that $\mathcal{P}(q)$ is, up to a fractional power of $q$ and an appropriate variable change, the multiplicative inverse of the Dedekind eta function:

\begin{align*}
\eta(\tau):=e^{\pi i\tau/12}\prod_{m=1}^{\infty}\left(1-e^{2\pi i m\tau}\right).
\end{align*}  Indeed, if we set $q:=e^{2\pi i\tau}$, then we have

\begin{align*}
\eta(\tau)^{-1} = \sum_{n=0}^{\infty}p(n)q^{n-1/24}.
\end{align*}  What makes $\eta(\tau)$ so remarkable is the following functional equation: for any $\gamma=\left(\begin{smallmatrix}
  a & b \\
  c & d 
 \end{smallmatrix}\right)\in\mathrm{SL}(2,\mathbb{Z})$,

\begin{align}
\eta\left(\frac{a\tau+b}{c\tau+d}\right) = \left( -i\left( c\tau+d \right) \right)^{1/2} \epsilon(a,b,c,d)\eta(\tau),\label{etamodularity}
\end{align} with

\begin{align*}
\epsilon(a,b,c,d) := \begin{cases} e^{b\pi i/12}, & c=0, d=1 \\ e^{\pi i\left( \frac{a+d}{12c} - s(d,c)\right)}, & c>0 \end{cases},
\end{align*} in which $s(h,k)$ is identical to the function in (\ref{pnformula}).  This functional symmetry is extremely unexpected, and difficult to prove.  It is the characteristic symmetry of the theory of modular forms, and it is this symmetry that allows for the construction of (\ref{pnformula}).

Certainly, these results raise even more questions to the arithmetician.  Given that $\left(p(n)\right)_{n\ge 0}$ is an integer sequence, it is very natural to ask when $p(n)$ is prime, or a square or higher power, or even or odd.  These sorts of questions cannot easily be answered with (\ref{pnformula}), which is not even clearly a convergent real number, to say nothing of its integrality.

What is truly astounding is the fact that the same theory of modular forms used to achieve (\ref{pnformula}) can be used---in a more refined way---to uncover \textit{arithmetic} information about $p(n)$.  In particular, Ramanujan discovered \cite{Ramanujan} the following properties (albeit with an error in his initial conjecture), which reveals a striking degree of structure to $p(n)$ despite its apparently pseudorandom appearance:

\begin{theorem}[Ramanujan, Watson, Atkin]\label{RamWA}
For all $\alpha\ge 1$ and $n\ge 0$, the following apply:
\begin{align}
\text{If } 24n\equiv 1\pmod{5^{\alpha}}\text{, then } p\left(n\right)&\equiv 0\pmod{5^{\alpha}};\label{RamWAp5}\\
\text{If } 24n\equiv 1\pmod{7^{\alpha}}\text{, then } p\left(n\right)&\equiv 0\pmod{7^{\left\lfloor\alpha/2\right\rfloor +1}};\label{RamWAp7}\\
\text{If } 24n\equiv 1\pmod{11^{\alpha}}\text{, then } p\left(n\right)&\equiv 0\pmod{11^{\alpha}}.\label{RamWAp11}
\end{align}
\end{theorem}

It is this result (and others like it) which we will focus on for the remainder of this paper.  Attempts to understand it have driven much of number theory in the twentieth century, and it continues to stimulate further research.  Notably, a true understanding of this theorem requires some understanding not only of arithmetic, but of complex analysis, abstract algebra, manifold theory, and topology.  Moreover, much of the theory is enormously intricate, and depends in large measure on computational and experimental work.

The remainder of this paper will be a description of how Theorem \ref{RamWA} is proved, how similar results have been discovered with respect to various more restrictive partition functions, the general difficulties in proving congruence families of this sort, and new approaches to the subject.  In Section 2 we give a sample of congruence families similar to those of Theorem \ref{RamWA}.  In so doing, we include some which are comparatively easy to prove, as well as some which cause difficulty of some sort to this day.

In Section 3 we will give a sketch of how to prove the simpler congruence families.  In particular, we focus on the first family of congruences in Theorem \ref{RamWA}.  This proof provides the primary template for proofs of other congruence families.  In Section 4 we will discuss many of the key complications that various specific families present to the proof method given in Section 3.  These complications include the genus and cusp count of the associated modular curve, representation of the associated modular function sequence in terms of eta quotients, and the existence of eigenfunctions modulo $\ell$ with respect to the associated linear operators $U_{\ell}$.  We point to the families given in Section 2 as examples of each complication.

We also wish to emphasize that this is a very interesting but small area of study in a much broader subject.  Questions about the arithmetic of $p(n)$ and associated functions have been pursued in substantially different directions by Ahlgren \cite{Ahlgren0}, Andrews \cite{AndrewsGarv}, Atkin \cite{AtkinL}, \cite{Atkin0}, Garvan \cite{AndrewsGarv}, \cite{Garvan0}, Lehner \cite{AtkinL}, Ono \cite{Ono2}, Radu \cite{Radu2}, and many others.

\section{Infinite Congruence Families}

Theorem \ref{RamWA} is by no means the only theorem of its kind.  Similar results apply to a very large variety of arithmetic sequences which are enumerated by a certain class of generating functions.  We will give an extremely brief---and far from comprehensive---list of results below, to give the reader a sense of the breadth of the arithmetic functions involved.

\subsection{The $j$ Invariant}\label{jinvariantsection}

We define the conventional Eisenstein series

\begin{align}
G_k(\tau) := \sum_{\substack{(m,n)\in\mathbb{Z}^2,\\ (m,n)\neq (0,0)}} \frac{1}{(m\tau+n)^{k}}.
\end{align}  These can easily be shown to be modular forms of weight $k$ for $k>2$.  We also define the modular discriminant $\Delta$ as the 24th power of $\eta(\tau)$ normalized with a power of $2\pi$:

\begin{align}
\Delta(\tau) := (2\pi)^{12}q\prod_{m=1}^{\infty}(1-q^m)^{24}.
\end{align}  With these functions, we may define the modular $j$ invariant:

\begin{align}
j:= j(\tau) = 1728\frac{60^3G_4(\tau)^3}{\Delta(\tau)} = \frac{1}{q}+744+\sum_{n=1}^{\infty}c(n)q^n.
\end{align}  See \cite[Chapter 1]{Diamond} for a standard treatment of this subject.

The $j$ invariant is one of the most important and foundational functions in the subject, and stands as a Hauptmodul for the entire modular group.  Lehner and Atkin have discovered an extraordinary range of arithmetic properties on its coefficients $c(n)$:

\begin{theorem}[Lehner, Atkin]
For $\alpha\ge 1$ and $n\ge 0$, the following apply:
\begin{align}
\text{If } n\equiv 0\pmod{2^{\alpha}}&\text{, then } c(n)\equiv 0\pmod{2^{3\alpha+8}};\label{jinvar2}\\
\text{If } n\equiv 0\pmod{3^{\alpha}}&\text{, then } c(n)\equiv 0\pmod{3^{2\alpha+3}};\label{jinvar3}\\
\text{If } n\equiv 0\pmod{5^{\alpha}}&\text{, then } c(n)\equiv 0\pmod{5^{\alpha+1}};\label{jinvar5}\\
\text{If } n\equiv 0\pmod{7^{\alpha}}&\text{, then } c(n)\equiv 0\pmod{7^{\alpha}};\label{jinvar7}\\
\text{If } n\equiv 0\pmod{11^{\alpha}}&\text{, then } c(n)\equiv 0\pmod{11^{\alpha}}\label{jinvar11}.
\end{align}
\end{theorem}  The families (\ref{jinvar2})-(\ref{jinvar7}) were proved by Lehner in \cite{Lehner0}, \cite{Lehner1}.  The more difficult family (\ref{jinvar11}) was proved by Atkin in \cite{Atkin}.  These properties are all the more interesting, given the importance of the coefficients $c(n)$ to applications in modular forms, group theory, Moonshine, and physics (e.g., \cite{Gannon}).

\subsection{$k$-Colored Partitions}\label{kcolorsection}

Many generalizations of $p(n)$ contain some analogous partition congruences.  For instance, we define the enumeration of partitions into $k$ colors as $p_{-k}(n)$, and give the generating function:

\begin{align}
P_{-k}(\tau) := \sum_{n=0}^{\infty} p_{-k}(n)q^n = \prod_{m=1}^{\infty}(1-q^{m})^{-k}.
\end{align}  In particular, we note that $P_{-1} = \mathcal{P}$, and that $p_{-1}(n)=p(n)$.  Dazhao Tang has noted a large number of interesting congruences for $p_{-k}(n)$ \cite{Tang}.  For example:

\begin{theorem}[D. Tang]
For $\alpha\ge 1$ and $n\ge 0$, the following applies:
\begin{align}
\text{If } 12n\equiv 1\pmod{5^{\alpha}}&\text{, then } p_{-2}(n)\equiv 0\pmod{5^{\left\lfloor\alpha/2\right\rfloor+1}}.
\end{align}
\end{theorem}

\subsection{Partitions Into Distinct Parts}\label{distinctpartssection}

Congruences similar in form to those of $p(n)$ often exist for more restrictive partition functions.  For example, if we define $p_D(n)$ as the counting function for partitions into distinct parts, then we of course have the generating function

\begin{align}
P_D(\tau) := \sum_{n=0}^{\infty} p_D(n)q^n = \prod_{m=1}^{\infty}(1+q^{m}) = \prod_{m=1}^{\infty}\frac{1-q^{2m}}{1-q^{m}}.
\end{align}  The following beautiful identity was proved by R\o dseth \cite{Rodseth} in 1969.

\begin{theorem}[R\o dseth]\label{distinctpartsthm}
For $\alpha\ge 1$ and $n\ge 0$, the following applies:
\begin{align}
\text{If } 24n\equiv -1\pmod{5^{2\alpha+1}}&\text{, then } p_{D}(n)\equiv 0\pmod{5^{\alpha}}.
\end{align}
\end{theorem}  \noindent Shane Chern and Mike Hirschhorn have since given a modified proof \cite{Chern}.

\subsection{Congruences Associated with Mock Theta Functions}\label{mockthetasection}

A particularly interesting contemporary topic is the quesiton of congruence families associated with mock theta functions.  Consider Ramanujan's third order mock theta function

\begin{align*}
\omega(q) := \sum_{n=0}^{\infty}\frac{q^{2n^2+2n}}{(q;q^2)_{n+1}^2}.
\end{align*}  Andrews, Dixit, and Yee have studied combinatorial interpretations of the coefficients of various mock theta functions \cite{AndDY}.  For example, they define $p_{\omega}(n)$ as the counting function of the number of partitions of $n$ in which all odd parts are less than twice the smallest part.  In that case, it can be shown that

\begin{align*}
\sum_{n=1}^{\infty} p_{\omega}(n)q^n = q\omega(q).
\end{align*}  As $p_{\omega}$ is a partition counting function, we can define a smallest-parts function for it, which we denote $\mathrm{spt}_{\omega}(n)$.  Liuquan Wang and Yifan Yang have proved \cite{Wang} that

\begin{theorem}[Wang, Yang]\label{wyfirst}
Let $\lambda_{\alpha}\in\mathbb{Z}$ be the minimal positive solution to $12x\equiv 1\pmod{5^{\alpha}}.$  Then
\begin{align}
\mathrm{spt}_{\omega}\left(2\cdot 5^{\alpha}n+\lambda_{\alpha}\right)\equiv 0\pmod{5^{\alpha}}.
\end{align}
\end{theorem}  This is proved by showing the family to be equivalent to another congruence family, this time corresponding to the generating function

\begin{align}
\sum_{n=0}^{\infty}c(n)q^n := \frac{2E_2(2\tau)-E_2(\tau)}{(q^2;q^2)_{\infty}},\label{E2}
\end{align} with $E_2(\tau)$ defined as the normalized holomorphic part of the weight 2 Eisenstein series:

\begin{align*}
E_2(\tau) := 1 - 24\sum_{n=1}^{\infty}\frac{nq^n}{1-q^n}.
\end{align*}  Theorem \ref{wyfirst} is in fact a corollary of the following:

\begin{theorem}\label{Thm12}
Let $12n\equiv 1\pmod{5^{\alpha}}$.  Then $c(n)\equiv 0\pmod{5^{\alpha}}$.
\end{theorem}  Wang and Yang gave the first proof.  A second proof using the localization method is given in \cite{Smoot0}.

\subsection{$k$-Elongated Plane Partitions}\label{kelongsection}

Another generalization of $p(n)$ is the $k$-elongated plane partition counting function, $d_k(n)$.  This was developed by George Andrews and Peter Paule as an application of their techniques for MacMahon's Partition Analysis program.

\begin{align}
D_k(\tau) := \sum_{n=0}^{\infty} d_k(n)q^n = \prod_{m=1}^{\infty}\frac{(1-q^{2m})^k}{(1-q^{m})^{3k+1}}.
\end{align}  Notice that $d_0(n)=p(n)$.

In 2021 George Andrews and Peter Paule conjectured \cite[Section 7.2, Conjecture 3]{AndrewsPaule} the following:

\begin{conjecture}[Andrews, Paule]\label{conjandpa}
For all integers $n\ge 0$, $\alpha\ge 1$, such that $8n\equiv 1\pmod{3^{\alpha}}$, we have
\begin{align}
d_2\left(n\right)\equiv 0\pmod{3^{\alpha}}.
\end{align}
\end{conjecture}  They suggested \cite[Section 8]{AndrewsPaule} that this conjecture would be difficult to prove.  In fact, the conjecture was refined and proved shortly thereafter \cite{Smoot2}, but only because the necessary methods have only very recently been developed (see Section \ref{cuspssection}).

\begin{theorem}[Smoot]\label{theoremsm}
For all integers $n\ge 0$, $\alpha\ge 1$, such that $8n\equiv 1\pmod{3^{\alpha}}$, we have
\begin{align}
d_2\left(n\right)\equiv 0\pmod{3^{2\left\lfloor\alpha/2\right\rfloor + 1}}.
\end{align}
\end{theorem}  A very large variety of congruence properties exist for $d_k(n)$, large even by the standards of partition functions, and a great deal of new results in this area are expected soon.  See, for example, the work of da Silva, Hirschhorn, and Sellers \cite{dasilvaet}.

\subsection{The Andrews--Sellers Congruences and Variations}\label{andsellvarsection}

A more ambitious problem lies in the study of the generalized 2-color Frobenius partition function $c\phi_2(n)$ first studied by Andrews in \cite{AndF}.  We do not give a combinatorial definition here, though we note that this is another generalization of $p(n)$.  We give the generating function

\begin{align}
\mathrm{C}\Phi_2(\tau) := \sum_{n=0}^{\infty} c\phi_2(n)q^n = \prod_{m=1}^{\infty}\frac{(1-q^{2m})^5}{(1-q^{m})^4(1-q^{4m})^2}.
\end{align}  James Sellers suggested the following congruence family in 1994 \cite{Sellers}, but it was not proved until Paule and Radu's work in 2012 \cite{Paule}:

\begin{theorem}[Paule, Radu]\label{cphialphathm}
For all integers $n\ge 0$, $\alpha\ge 1$, such that $12n\equiv 1\pmod{5^{\alpha}}$, we have
\begin{align}
c\phi_2\left(n\right)\equiv 0\pmod{5^{\alpha}}.\label{defnLalpha}
\end{align}
\end{theorem}  This was an extremely difficult result to prove, and required the development of new techniques which we will discuss later.  This was all the more stunning given the fact that it was proposed some eighty years after Ramanujan's original conjecture.

\section{Basic Proof Method}

As we have remarked, there exists a striking variety in the examples above in terms of the difficulty of proofs.  Nevertheless, we can recognize a foundational approach to all of these examples.  This approach succeeds in proving some of the congruence families above.  For the more resistant examples, we can modify this approach by accounting for various complications which we give in the sequel.

\subsection{Linear Congruences for $p(n)$}

Euler examined $p(n)$ as early as 1750 without discovering any obvious arithmetic properties.  This makes Ramanujan's work in this subject all the more remarkable.  We will begin with the initial cases of Theorem \ref{RamWA}, which Ramanujan discovered upon studying tables of $p(n)$:

\begin{theorem}[Ramanujan]\label{linearcases}
Let $n\in\mathbb{Z}_{\ge 0}$.  Then
\begin{align}
p\left(5n+4\right)&\equiv 0\pmod{5}\\
p\left(7n+5\right)&\equiv 0\pmod{7}\\
p\left(11n+6\right)&\equiv 0\pmod{11}.
\end{align}
\end{theorem}  In retrospect, these relations---the first in particular---are especially easy to check, and it is striking that earlier mathematicians missed them.  We will take the first relation as a useful example.  There are multiple ways of proving the congruence, but we are especially interested in using the following identity:

\begin{theorem}[Ramanujan]
Let $q\in\mathbb{C}$ such that $|q|<1$.  Then
\begin{align}
\sum_{n=0}^{\infty}p(5n+4)q^n = 5\prod_{m=1}^{\infty}\frac{(1-q^{5m})^5}{(1-q^m)^6} = 5\frac{\mathcal{P}(q)^6}{\mathcal{P}(q^5)^5}.\label{RKor}
\end{align}
\end{theorem}  As with $\mathcal{P}$, we may easily verify the existence of an integer power series representation of the right-hand side by a geometric series expansion of each factor.  Determining convergence is not difficult, although this identity can also be proved in terms of $q$ as a formal indeterminate.  We will provide a sketch of an analytic proof which reveals the first instance of modular forms playing a role in determining the arithmetic of $p(n)$.  Indeed, one can ``normalize" the identity to read

\begin{align}
\prod_{m=1}^{\infty}(1-q^{5m})\cdot \sum_{n=0}^{\infty}p(5n+4)q^{n+1} = 5q\prod_{m=1}^{\infty}\left(\frac{1-q^{5m}}{1-q^m}\right)^6.\label{RKorA}
\end{align}  We denote the right-hand side (without the factor 5) as

\begin{align}
t:=t(\tau) = q\prod_{m=1}^{\infty}\left(\frac{1-q^{5m}}{1-q^m}\right)^6,
\end{align} with $q=e^{2\pi i\tau}$, and $\tau\in\mathbb{H}$.  Notice that this domain restriction for $\tau$ implies that $|q|<1$.  Also notice that, using our definition of $\eta$, we have

\begin{align}
t = \left(\frac{\eta(5\tau)}{\eta(\tau)}\right)^6.
\end{align}  This function posesses a much cleaner form of symmetry than that shown by $\eta(\tau)$ alone.  Indeed, for any $N\in\mathbb{N}$, let us define the \textit{congruence subgroup}

\begin{align*}
\Gamma_0(N) = \Bigg\{ \begin{pmatrix}
  a & b \\
  c & d 
 \end{pmatrix}\in \mathrm{SL}(2,\mathbb{Z}) : N|c \Bigg\}.
\end{align*}  We also define a group action $\Gamma_0(N)\times\mathbb{H}\rightarrow\mathbb{H}$ by

\begin{align*}
\gamma\tau := \frac{a\tau+b}{c\tau+d}
\end{align*} for $\gamma=\left(\begin{smallmatrix}
  a & b \\
  c & d 
 \end{smallmatrix}\right)\in\Gamma_0(N)$, $\tau\in\mathbb{H}$.

In that case, for $\gamma\in\Gamma_0(5)$, we have

\begin{align}
t\left(\gamma\tau\right) = t(\tau).
\end{align}

On the other hand, let us consider $\rho := e^{2\pi i/5}$, and define $q^{1/5} := e^{2\pi i\tau/5}$, and $\kappa\in\mathbb{Z}$.  In that case we have

\begin{align}
&\sum_{\lambda=0}^4 \rho^{-4\kappa\lambda} q^{-4/5} \mathcal{P}\left( \rho^{\kappa\lambda}q^{1/5} \right) \label{p4np5weightsum}\\
&= \sum_{\lambda=0}^4 \sum_{n=0}^{\infty} e^{-2\pi i\cdot 4\kappa\lambda/5} e^{-2\pi i\tau\cdot 4/5} p(n) e^{2\pi i n\left(\tau+\kappa\lambda\right)/5}\\
&= \sum_{n=0}^{\infty}p(n) e^{2\pi i \tau (n-4)/5 } \sum_{\lambda=0}^4  e^{2\pi i \kappa\lambda (n-4)/5}\\
&= 5\sum_{n=0}^{\infty}p(5n+4) q^n.
\end{align}  This gives us a precise expression of $\sum_{n=0}^{\infty}p(5n+4) q^n$ that can be suitably manipulated.  Indeed, rewriting (\ref{p4np5weightsum}) in terms of eta functions, we have

\begin{align*}
\sum_{\lambda=0}^4 \rho^{-4\kappa\lambda} q^{-4/5} \mathcal{P}\left( \rho^{\kappa\lambda}q^{1/5} \right) &= \sum_{\lambda=0}^4 \rho^{\lambda\kappa(-24(4)+1)/24}q^{(1-24(4))/(24\cdot 5)}\eta\left( \frac{\tau+\kappa\lambda}{5} \right)^{-1},
\end{align*} so that

\begin{align*}
q^{(24(4)-1)/(24\cdot 5)}\sum_{n=0}^{\infty}p(5n+4) q^n &= \frac{1}{5} \sum_{\lambda=0}^4 \rho^{\lambda\kappa(5)/24}\eta\left( \frac{\tau+\kappa\lambda}{5} \right)^{-1}.
\end{align*}  If we take $\kappa=24$, and simplify, then we have

\begin{align*}
h(\tau) := q^{19/24}\sum_{n=0}^{\infty}p(5n+4) q^n &= \frac{1}{5} \sum_{\lambda=0}^4 \rho^{5\lambda}\eta\left( \frac{\tau+24\lambda}{5} \right)^{-1}.
\end{align*}  Our understanding of the modular symmetry of $\eta$ suggests that $h(\tau)$ posesses a similar symmetry.  Indeed, if we multiply by $\eta(5\tau)$ and define

\begin{align*}
L_1(\tau):=\eta(5\tau)\cdot h(\tau) = \prod_{m=1}^{\infty}\left( 1-q^{5m} \right)\cdot \sum_{n=0}^{\infty}p(5n+4) q^{n+1} &= \frac{1}{5} \sum_{\lambda=0}^4 \rho^{5\lambda}\frac{\eta(5\tau)}{\eta\left( \frac{\tau+24\lambda}{5} \right)},
\end{align*} then it can be shown that $L_1(\gamma\tau)=L_1(\tau)$ for $\gamma\in\Gamma_0(5)$.

The significance of the symmetric behavior of $t$ over the group action of a subgroup of $\mathrm{SL}(2,\mathbb{Z})$ arises from multiple important points.

Much of the theory of complex manifolds has arisen from the study of such transformations on $\mathbb{H}$.  For a fixed $N$, suppose we take the quotient of $\mathbb{H}$ by the group action of $\Gamma_0(N)$.  That is, we consider any two points $\tau_1,\tau_2$ to be equivalent if there exists some $\gamma\in\Gamma_0(N)$ such that $\tau_1=\gamma\tau_2$.

The resulting manifolds can be compactified by the addition of a finite number of additional points, which we will call \textit{cusps}.  The idea is that we also identify points of $\mathbb{Q}$ with respect to the group action.  This partitions all members of $\mathbb{Q}$ into a finite number of subsets.  In particular, the subset containing fractions identified with $0/1$ is called the cusp at 0, and is usually denoted $[0]_N$.

In addition, if we extend $\mathbb{Q}$ to include expressions $k/0$ for $k\neq 0$, then we find such expressions identified with fractions in which the denominators are divisible by $N$.  This subset is called the cusp at $\infty$, and is denoted $[\infty]_N$.

Taking into account the cusp $[\infty]_N$, together with the finite number of remaining cusps, we have a compact manifold which is referred to as the classical modular curve, $\mathrm{X}_0(N)$.  In the case that $N$ is a prime number (like 5), the associated curve is compactified by the addition of only two cusps: $[0]_N$ and $[\infty]_N$.

Notice that the functions $t, L_1$ each induce certain well-defined functions $\hat{t}$, resp. $\hat{L_1}$, on $\mathrm{X}_0(5)$.  Moreover, both of these induced functions posesses an analogue of \textit{meromorphicity} on $\mathrm{X}_0(5)$.  In particular, they are holomorphic along all of $\mathrm{X}_0(5)$ except for the cusp $[0]_5$, at which they both posess the analogue of a principal part.

All of this may seem somewhat ostentatious; however, this subject provides us with an extremely important theorem, which dictates that \textit{no nonconstant function can be holomorphic everywhere on a compact Riemann surface.}  This gives rise to certain finiteness conditions which are extremely powerful when it comes to verifying  identities between such functions.

In our case, $5\hat{t}$, $\hat{L_1}$ have matching principal parts at $[0]_5$ of $\mathrm{X}_0(5)$.  Therefore, the function $\hat{L_1}-5\cdot\hat{t}$ is holomorphic over the entirety of $\mathrm{X}_0(5)$, and is therefore a constant.  It can be quickly verified that the constant term in a power series expansion of $\hat{L_1}-5\hat{t}$ is 0.  Thus, $\hat{L_1}-5\hat{t}=0$.  However, $\hat{L_1}-5\hat{t}$ is induced from $L_1-5t$.  So if $\hat{L_1}-5\hat{t}= 0$, then we must also have $L_1-5t\neq 0$.  This gives us (\ref{RKorA}), and (\ref{RKor}).

The functions $t, L_1$ each constitute \textit{modular functions over the subgroup} $\Gamma_0(5)$.  More generally, a modular function $f$ over $\Gamma_0(N)$ is a function holomorphic over $\mathbb{H}$ which exhibits $f(\gamma\tau)=f(\tau)$ for all $\gamma\in\Gamma_0(N)$ and $\tau\in\mathbb{H}$, and whose induced function $\hat{f}$ on $\mathrm{X}_0(N)$ exhibits meromorphic behavior at the cusps.

If we define $\mathcal{M}\left(\Gamma_0(N)\right)$ to be the set of modular functions over $\Gamma_0(N)$, then (\ref{RKor}) can be understood as two different representations of a member of $\mathcal{M}\left(\Gamma_0(5)\right)$.

Similar arguments can be used to justify the congruences modulo 7, 11 in Theorem \ref{linearcases}, although the case modulo 11 is much more difficult.

\subsection{Infinite Congruence Families}

Other techniques besides that outlined above can be used to verify the results of Theorem \ref{linearcases}.  Nevertheless, there are multiple reasons for preferring the modular approach.  In the first case, the theory is so well-understood that it lends itself well to generalization and algorithmization.  However, the most compelling reason is that this approach can be used to reveal far richer arithmetic structure.

After discovering the three linear cases above, Ramanujan also discovered that

\begin{align}
p\left(25n+24\right)&\equiv 0\pmod{25}\\
p\left(49n+47\right)&\equiv 0\pmod{49},
\end{align} and he found evidence that
\begin{align}
p\left(121n+116\right)&\equiv 0\pmod{121}.
\end{align}  These results led him to suggest the general infinite family which eventually became Theorem \ref{RamWA}.

These congruence families represent especially deep results.  Let us focus for the time being on the family (\ref{RamWAp5}).

\begin{align}
p\left(5^{\alpha}n+\lambda_{5,\alpha}\right)&\equiv 0\pmod{5^{\alpha}}.
\end{align}  Notice that we proved the first case of this congruence family by constructing a modified generating function for $p\left(5n+\lambda_{5,1}\right)$, which we denoted as $L_1$.  We then showed that $L_1$ is a modular function, and found a representation of $L_1$ in which the divisibility by 5 was apparent.  This approach is especially useful in approaching the general family of congruences.  For $\alpha\ge 1$ we define

\begin{align}
L_{\alpha} := \phi_{\alpha}(q)\cdot\sum_{n=0}^{\infty}p\left(5^{\alpha}n+\lambda_{5,\alpha}\right)q^{n+1},
\end{align} where

\begin{align}
\phi_{2\alpha-1} :=& \prod_{m=1}^{\infty}\left( 1-q^{5m} \right),\\
\phi_{2\alpha} :=& \prod_{m=1}^{\infty}\left( 1-q^{m} \right).
\end{align}  By similar steps to those in the preceeding section, we may show that $L_{\alpha}$ is modular over $\Gamma_0(5)$.  An inductive approach naturally suggests itself.  Let us also define

\begin{align}
\mathcal{A} := \frac{\eta(25\tau)}{\eta(\tau)} =& q\prod_{m=1}^{\infty}\frac{1-q^{25m}}{1-q^m} = q\frac{\mathcal{P}(q)}{\mathcal{P}(q^{25})}.
\end{align}  We then define the following linear operators:

\begin{align}
U_5^{(1)}\left(\sum_{n\ge N}a(n)q^n\right) &:= \sum_{5n\ge N}a(5n)q^n,\\
U_5^{(0)}\left(\sum_{n\ge N}a(n)q^n\right) &:= U_5^{(1)}\left(\mathcal{A}\cdot\sum_{n\ge N}a(n)q^n\right),\\
U_5^{(\alpha)}\left(\sum_{n\ge N}a(n)q^n\right) &:= U_5^{(\alpha\bmod{2})}\left(\sum_{n\ge N}a(n)q^n\right).
\end{align}

This construction is straightforward to the theory, though we briefly cover it for want of space.  One can then show that

\begin{align}
U_5^{(\alpha)}\left( L_{\alpha} \right) = L_{\alpha+1}.
\end{align}  From here, the proof method is clear: we need to show that if $5^{\alpha}|L_{\alpha}$, then upon applying $U_5^{(\alpha)}$, we have $5^{\alpha+1}|L_{\alpha+1}$.

To do this, we need to find a suitable space $\mathcal{K}$ that is general enough that $L_{\alpha}\in\mathcal{K}$, yet specific enough that we can extract interesting arithmetic information.  In the case of powers of 5 in Theorem \ref{RamWA}, we found that

\begin{align*}
L_{1}=5t.
\end{align*}  More generally, it turns out that for all $\alpha\ge 1$,

\begin{align}
L_{\alpha}\in\mathbb{Z}[t].
\end{align}  The curious role played by $t$ is explained by the fact that the space of all meromorphic functions over $\mathrm{X}_0(5)$ with a pole only at the cusp $[0]_5$ is equal to the polynomial ring over $\mathbb{C}[\ \hat{t}\ ]$.  As a result, we have

\begin{align}
\mathcal{M}^{0}\left( \Gamma_0(5) \right) = \mathbb{C}[t].
\end{align}  The induced functions $\hat{L}_{\alpha}$ on $\mathrm{X}_0(5)$ have poles only at $[0]_5$.

\subsection{Completing the Proof}

The importance of determining the algebraic structure of $\mathcal{M}^{0}\left( \Gamma_0(5) \right)$ is that our problem is reduced in this case to a question of how $U_5^{(\alpha)}(t^n)$ may be expressed.  It is at this point that we may use an important result for $t$ to build a recurrence.

\begin{theorem}
Define
\begin{align*}
a_4(\tau) &=-(5^{12}t^5+6\cdot 5^{10}t^4+63\cdot 5^{7}t^3+52\cdot 5^5 t^4+63\cdot 5^2 t^5)\\
a_3(\tau) &= -(5^{9}t^4+6\cdot 5^{7}t^3+63\cdot 5^{4}t^2+52\cdot 5^{2}t)\\
a_2(\tau) &=-(5^6t^3+6\cdot 5^4 t^2+63\cdot 5 t)\\
a_1(\tau) &= -(5^3 t^3+6\cdot 5 t)\\
a_0(\tau) &= -t.
\end{align*}  

\begin{align}
t^5+\sum_{j=0}^4 a_j(5\tau) t^j = 0.\label{modX}
\end{align}
\end{theorem}  If we add to this the fact that that $U_5(f(5\tau)\cdot g(\tau)) = f(\tau)\cdot U_5(g(\tau))$, then we have

\begin{align}
U_5^{(\alpha)}\left(t^m\right) = -\sum_{j=0}^4 a_j(\tau) U_5^{(\alpha)}\left(t^{m+j-5}\right).\label{modXa}
\end{align}  As such, we need only directly compute $U_5^{(\alpha)}(t^m)$ for only a finite number of $m$.

Moreover, we can use (\ref{modXa}) to construct lower bounds for the 5-adic valuations in the polynomial expansions of $U_5^{(\alpha)}\left(t^m\right)$.  One can show that, for certain auxiliary functions $h_i(m,r)$ with finite support,

\begin{align}
U_5^{(0)}\left(t^m\right) &= \sum_{r\ge\left\lceil (m+1)/5 \right\rceil} h_0(m,r) 5^{\left\lfloor\frac{5r-m-2}{2}\right\rfloor} t^r\label{u50tn}\\
U_5^{(1)}\left(t^m\right) &= \sum_{r\ge\left\lceil m/5 \right\rceil} h_1(m,r) 5^{\left\lfloor\frac{5r-m-1}{2}\right\rfloor} t^r\label{u51tn}.
\end{align}

From here, we need to apply the operators $U_5^{(\alpha)}$ on the functions $L_{\alpha}$.  The resulting function will be a polynomial in $t$, whose coefficients will be divisible by a nonnegative power of 5.  We need to keep a sufficiently high lower bound on the 5-adic valuation of these coefficients.  

For $i=0,1$ we define the functions

\begin{align*}
\theta_i(m) &:= \left\lfloor \frac{5m-i}{2} \right\rfloor,
\end{align*} and the polynomial spaces

\begin{align}
\mathcal{V}^{(i)} :=& \left\{ \sum_{m\ge 1} s(m)\cdot 5^{\theta_i(m)}\cdot t^m : s(m)\text{ has finite support} \right\}.
\end{align}  Now let us take some $f\in\mathcal{V}^{(0)}$.  Then we have

\begin{align}
f=\sum_{m\ge 1} s(m)\cdot 5^{\theta_0(m)}\cdot t^m.
\end{align}  If we apply $U_5^{(0)}$ to $f$ and expand, then

\begin{align}
U_5^{(0)}(f)&=\sum_{m\ge 1} s(m)\cdot 5^{\theta_0(m)}\cdot U_5^{(0)}(t^m)\\
&=\sum_{m\ge 1}\sum_{r\ge\left\lceil (m+1)/5 \right\rceil} s(m)h_0(m,r)\cdot 5^{\theta_0(m)+\left\lfloor\frac{5r-m-2}{2}\right\rfloor} t^r.
\end{align}  If we examine the power of 5, we find that

\begin{align}
\theta_0(m)+\left\lfloor\frac{5r-m-2}{2}\right\rfloor &= \left\lfloor \frac{5m}{2} \right\rfloor + \left\lfloor\frac{5r-m-2}{2}\right\rfloor\\
&\ge \left\lfloor\frac{5r+4m-3}{2}\right\rfloor\\
&\ge \left\lfloor\frac{5r-1}{2}\right\rfloor + \left\lfloor\frac{4m-2}{2}\right\rfloor\\
&\ge \theta_1(r) + 1.
\end{align}  Therefore, we must have

\begin{align}
\frac{1}{5}U_5^{(0)}(f)\in \mathcal{V}^{(1)}.
\end{align}  Similarly, we can prove that if $f\in\mathcal{V}^{(1)}$, then

\begin{align}
\frac{1}{5}U_5^{(1)}(f)\in \mathcal{V}^{(0)}.
\end{align}

Finally, we note that 

\begin{align}
\frac{1}{5}L_1 = t\in \mathcal{V}^{(1)},
\end{align} and that $L_2 = U_5^{(1)}\left(L_1\right)$.  This is enough to complete the proof.  The sequence $\left( L_{\alpha} \right)_{\alpha\ge 1}$ alternates between $\mathcal{V}^{(1)}$ and $\mathcal{V}^{(0)}$, gaining a power of 5 upon each application of $U_5^{(\alpha)}$.

For additional details, we direct the reader to \cite{Paule2} for a modern treatment.  For a classical treatment, we refer the reader to Watson's original paper \cite{Watson} (in German), or \cite[Chapters 7, 8]{Knopp}.

\section{Key Difficulties}

Our focus on (\ref{RamWAp5}) is justified by the fact that this is in many ways the most straightforward congruence family to undrerstand.  Watson published the first proof in 1938 \cite{Watson}, but Ramanujan himself appears to have understood the proof technique \cite{BO}.

The second part of Theorem \ref{RamWA}, (\ref{RamWAp7}), was only slightly more difficult to prove.  Ramanujan actually believed that the congruence family had identical form to that of powers of 5.  Because he was not aware of the half-growth that we see in the complete theorem, we generally give Watson credit for being the first person to understand and complete the proof in 1938 \cite{Watson}.

The third part of Theorem \ref{RamWA}, (\ref{RamWAp11}), is much, much more difficult.  Certainly, Ramanujan did not have a proof, and no proof was found until Atkin's work in 1967 \cite{Atkin}.  This is almost closer to our time than to Ramanujan's.

As we showed in Section 2, a large number of similar congruence families turn out to exist for a great variety of different arithmetic functions.  What is especially strange is the enormous variety with respect to the \textit{difficulty} of the associated proofs.  Given that the outline above was understood by at least one person more than a century ago, it is striking that many other congruence families have taken decades to prove, and even some contemporary families of congruences are extremely difficult to prove.  On the other hand, some congruence families are so well-understood that they can be proved routinely in mathematics journals by the methods described in the previous section.

We want to discuss what the possible complications are.

\subsection{Genus}

The first major complication, and in some ways the easiest to understand, is the problem involving the \textit{genus} of the associated modular curve.  The genus of the curve $\mathrm{X}_0(N)$ is denoted $\mathfrak{g}\left( \mathrm{X}_0(N) \right)$.  This is of course a straightforward topological number: namely, the number of handles in the compact manifold.  A sphere has genus 0, while a torus or coffee mug has genus 1.

It is particulary strange that the topology of $\mathrm{X}_0(N)$ determines the difficulty in proving an associated partition congruence family.  However, as in the case of proving (\ref{RamWAp5}) from Theorem \ref{RamWA}, we express an interest in proving a given congruence family by expressing the associated functions $L_{\alpha}$ in terms of some reference functions.  The most convenient reference functions are generally the generators of the space of modular functions whose induced equivalent functions are holomorphic across the entire associated modular curve, except at a single point.  Such a space turns out to have a very close relationship with $\mathfrak{g}$.

We already discussed the extremely important theorem that a nonconstant meromorphic function on $\mathrm{X}_0(N)$ must have at least one pole somehwere on the surface.  Let us consider the set of all functions $\hat{\mathcal{M}}^p\left( \mathrm{X}_0(N) \right)$ which are holomorphic on $\mathrm{X}_0(N)$ except for a single pole at the point $[p]_N$.  Take some $f\in\hat{\mathcal{M}}^p\left( \mathrm{X}_0(N) \right)$.  We ask what possible order can $f$ take at $[p]_N$.

\begin{theorem}[Weierstrass]\label{WGT}
Let $\mathrm{X}$ be a compact Riemann surface, and let 
\begin{align*}
f:\mathrm{X}\longrightarrow\hat{\mathbb{C}}
\end{align*} be holomorphic over $\mathrm{X}$, except for a pole at a point $p\in\mathrm{X}$.  Then the order of $f$ at $p$ can assume any negative integer, with exactly $\mathfrak{g}\left( \mathrm{X} \right)$ exceptions, which must be members of the set $\{1, 2, ..., 2\mathfrak{g}-1\}$.
\end{theorem}

We see that the genus plays a direct role in the complexity of the space of functions that we will end up working over.  In particular, if $\mathfrak{g}\left( \mathrm{X}_0(N) \right)=0$, then we will have

\begin{align}
\mathcal{M}^{p}\left( \Gamma_0(N) \right) = \mathbb{C}[x],
\end{align} for some function $x$ (here called a Hauptmodul).  On the other hand, if $\mathfrak{g}\left( \mathrm{X}_0(N) \right)=1$, then there must be exactly a single exception to the possible orders of $f$ at $p$.  But of course, if $f$ has order $-1$ at $p$, then $f^n$ will have order $-n$ at $p$ for any $-n\in\mathbb{Z}_{\le -1}$ that we like, and no exception would exist.  Therefore, the exceptional order must be $-1$.  As a result, the simplest algebraic representation of the space of functions at $p$ must have the form

\begin{align}
\mathcal{M}^{p}\left( \Gamma_0(N) \right) = \mathbb{C}[x]\oplus y\cdot\mathbb{C}[x],
\end{align} in which $x$ has order $-2$ and $y$ has order $-3$.

This immediately complicates any proof that we intend to construct for the given congruence family.  If we wanted to examine the representation of $L_{\alpha}\in\mathbb{C}[x]\oplus y\cdot\mathbb{C}[x]$, then we would need to examine $U_{\ell}^{(\alpha)}(x^n)$, $U_{\ell}^{(\alpha)}(y\cdot x^n)$.  This generally requires that we examine twice as many specific relations for a finite range of $n$ as we would expect for a genus 0 surface.  Moreover, this also poses complications for a given proof, since $U_{\ell}^{(\alpha)}(y^m\cdot x^n)$ will have to be expressed in terms of $y^jx^k$.

This is one major reason that (\ref{RamWAp11}) in Theorem \ref{RamWA} is so difficult: the associated modular curve $\mathrm{X}_0(11)$ has genus 1, in contrast to the modular curves associated with (\ref{RamWAp5}) and (\ref{RamWAp7}), each of which is associated with a genus 0 curve.

The same problem applies to the family (\ref{jinvar11}) for the coefficients of the $j$-invariant, although (\ref{jinvar2})-(\ref{jinvar7}) are much easier to prove.

A third extraordinarily difficult example lies in the Andrews--Sellers congruence family, which is associated with the genus 1 curve $\mathrm{X}_0(20)$.

These examples are both plagued with a second major difficulty.  The former two examples suffer from a lack of a suitable eta quotient basis, while the latter suffers from a much deeper and more challenging problem.  We will briefly discuss both below.

\subsection{Number of Cusps}\label{cuspssection}

A second key complication which has only recently become studied, is the number of cusps of $\mathrm{X}_0(N)$, denoted $\epsilon\left( \mathrm{X}_0(N) \right)$.  In particular, if $N$ is prime, then $\epsilon\left( \mathrm{X}_0(N) \right)=2$.  However, for composite $N$, $\epsilon\left( \mathrm{X}_0(N) \right)$ will generally be larger than 2.  For example, $\epsilon\left( \mathrm{X}_0(10) \right)=4$.

The importance of this complication emerges in the question of the function sequence $\left( L_{\alpha} \right)_{\alpha\ge 1}$ associated with the congruence family.  The induced function $\hat{L_{\alpha}}$ will only have poles at the cusps.  In the case that there are only 2 cusps, $\hat{L_{\alpha}}$ must have a pole at one cusp and a zero at the other cusp.  This greatly simplifies the question of representing $\hat{L_{\alpha}}$, and therefore $L_{\alpha}$.  On the other hand, if there are more than 2 cusps, then $\hat{L_{\alpha}}$ may exhibit more complexity in terms of where its poles are.

The most straightforward solution to this problem is to take advantage of the fact that we can always push a given modular function $f$ into the space of functions which are holomorphic everywhere except at a single cusp (say, $[0]_N$), by multiplying by a certain computable modular function $z$.  In the case that the associated modular curve $\mathrm{X}_0(N)$ has genus 0, we have

\begin{align}
z,\ z\cdot f\in\mathcal{M}^{0}\left( \Gamma_0(N) \right) = \mathbb{C}[x],
\end{align} for a suitable Hauptmodul $x$.  Dividing by $z$, we have

\begin{align}
f&\in\mathbb{C}[x]_{\mathcal{S}},
\end{align} in which $\mathbb{C}[x]_{\mathcal{S}}$ is the localization of the polynomial ring at the multiplicatively closed set
\begin{align}
\mathcal{S}&=\{z^n : n\ge 0\}.
\end{align}  This gives us an expression of $f$ as a rational polynomial in terms of $x$.  With the right choice of $x$ and $z$, this can exhibit the key arithmetic information needed.

As an example, we consider Conjecture \ref{conjandpa} regarding the 2-elongated plane partition function $d_2(n)$.  This congruence family is associated with the modular curve $\mathrm{X}_0(6)$.  This curve has genus 0, but it has four distinct cusps.  Indeed, for the associated function sequence $\left( L_{\alpha} \right)_{\alpha\ge 1}$, we have

\begin{align}
L_{\alpha}\in\mathcal{M}\left( \Gamma_0(6) \right)\setminus\mathcal{M}^{0}\left( \Gamma_0(6) \right).
\end{align}  With the appropriately chosen eta quotient $x$, we have

\begin{align}
L_{\alpha}\in\mathbb{Z}[x]_{\mathcal{S}},
\end{align} with

\begin{align}
\mathcal{S}:=\left\{(1+9x)^n : n\ge 0\right\},
\end{align} for all $\alpha\ge 1$.  This representation of $L_{\alpha}$ yields the divisibility properties inherent in the congruence family.  Once the cusp count is accounted for, the proof is straightforward (see \cite{Smoot2}).

This technique is denoted as the localization method for proving congruence families.  It has also been used to give an alternative proof of the congruences in Theorem \ref{wyfirst} \cite{Smoot0}.  On the other hand, Theorem \ref{distinctpartsthm} was proved by other techniques (see Section \ref{eigenfunctionsection}), although a proof by localization ought to be possible.

The appeal of the localization approach is that it expresses the associated congruence family in a very natural environment, usually involving reference functions which can be derived by a straightforward process.  Indeed, in many cases the theory is so well understood that the prospect of an algorithmic approach to the problem appears promising.  The disadvantage is that the steps in such a proof are often lengthy, although this would pose little problem if an algorithmic procedure is one day developed using these techniques.

\subsection{Existence of an Eta Quotient Basis}

The two complications previously discussed are directly tied to the topology of the modular curve associated with the given congruence family.  The next two complications are less obviously tied to the topology of the modular curve; nevertheless, they are exceedingly difficult complications that need to be addressed.

We already referred to the case for powers of 11 in Theorem \ref{RamWA}.  The associated modular curve $\mathrm{X}_0(11)$ has genus 1, which immediately complicates the relevant space of modular functions.  On the other hand, because 11 is prime, only two cusps exist.  As such, we know that the associated function sequence $\left( L_{\alpha} \right)_{\alpha\ge 1}$ will exist in the ring

\begin{align}
\mathcal{M}^{0}\left( \Gamma_0(11) \right) = \mathbb{C}[x]\oplus y\mathbb{C}[x],\label{possiblerank211}
\end{align} in which $x,y$ are functions with orders $-2, -3$, respectively, at the cusp $[0]$.  However, if we examine the space of eta quotients over $\Gamma_0(11)$, we find that no such eta quotient exists.  The functions $x,y$ certainly exist; but they do not admit such simple representations.  As an example, we consider the first case of the third congruence family (i.e., powers of 11) in Theorem \ref{RamWA}.  If we attempt to find an analogue of Theorem \ref{RKor} for $p(11n+6)$, we discover a much more cumbersome identity which is not easy to express in terms of reference functions.  One example, taken from \cite{Paule2} and closely related to Atkin's work in \cite{Atkin}, is the following:

\begin{align}
\prod_{m=1}^{\infty}&(1-q^{11m})\sum_{n=0}^{\infty}p(11n+6)q^{n+1}\nonumber\\ = & 11^3\cdot\frac{1}{2}\cdot \left( 2\cdot f_3^2 + U_2\left(f_2^2\right) \right)+ 11\cdot\frac{7}{2}\left( f_2 - 4\cdot U_2\left( f_3 \right) \right) - 11\left( f_3 - U_2\left( f_2 \right) \right) + 11^4 z,\label{pn116succ}
\end{align} with

\begin{align*}
f_2 &:= q^2 \prod_{m=1}^{\infty}\frac{(1-q^{11n})^3(1-q^{22n})}{(1-q^n)(1-q^{2n})^3},\\
f_3 &:= q^3 \prod_{m=1}^{\infty}\frac{(1-q^{11n})(1-q^{22n})^3}{(1-q^n)^3(1-q^{2n})},\\
z &:= q^5 \prod_{m=1}^{\infty}\left(\frac{1-q^{11m}}{1-q^m}\right)^{12}.
\end{align*}

The functions $f_2, f_3$ are not modular for $\Gamma_0(11)$, but rather for $\Gamma_0(22)$.  To properly form the functions $x,y$ corresponding to (\ref{possiblerank211}), we need to define

\begin{align*}
x &:= U_2(f_2) - f_3,\\
y &:= \frac{1}{11}\left( 4 U_2(f_3) - f_2 \right).
\end{align*}  Notice that the presence of the $U_2$ operator suggests (correctly) that $x,y$ cannot simply be expressed in terms of $\eta(\tau)$, $\eta(11\tau)$ alone.  This adds to the complications of proving Ramanujan's congruence family for powers of 11.  Indeed, the first proof of this congruence family was done by Atkin, who chose a much more complex ring structure for the associated modular functions.

We note also that this complication also plagues the congruences (\ref{jinvar11}) for the coefficients of the $j$-invariant.

\subsection{Existence of Eigenfunctions Modulo $\ell$}\label{eigenfunctionsection}

We now consider an especially difficult complication in the subject of congruence families---indeed, perhaps the most difficult complication that we know of.

Let us consider the case of the Andrews--Sellers congruences, in which the associated sequence of modular functions exists over $\Gamma_0(20)$.  We have 

\begin{align*}
L_0 &= 1,\\
L_{\alpha} &= \Phi_{\alpha}\cdot \sum_{12n\equiv 1\bmod{5^{\alpha}}} c\phi_2(n)q^{\left\lfloor n/5^{\alpha} \right\rfloor},
\end{align*} with $\Phi_{\alpha}$ suitably chosen.  To prove that $a(n)\equiv 0\pmod{5^{\alpha}}$ whenever $12n\equiv 1\pmod{5^{\alpha}}$, we need to argue that for every $M\in\mathbb{Z}_{>0}$ there exists an $N\in\mathbb{Z}_{>0}$ such that for all $n\ge N$,

\begin{align*}
L_n\equiv 0\pmod{5^M}.
\end{align*}  In particular, we need to show that $N=M$ will suffice.

Because the genus of $\mathrm{X}_0(20)$ is 1, we might ordinarily expect to simply express each function $L_{\alpha}$ as a member of the rank 2 $\mathbb{C}[x]$-module

\begin{align}
\mathcal{M}^{0}\left( \Gamma_0(20) \right) = \mathbb{C}[x]\oplus y\cdot\mathbb{C}[x],
\end{align} with integer coefficients.  One needs only to consider how the associated $U^{(\alpha)}$ operators act on $x^n$, $yx^n$.  However, at this point a new and extraordinary complication emerges.

To understand this complication, we define the functions $x$, $y$ which give us the rank 2 module above as the following:

\begin{align}
x = x(\tau) :=& q\prod_{m=1}^{\infty}\frac{(1-q^{2m})(1-q^{10m})^3}{(1-q^m)^3(1-q^{5m})}\label{ydef}\\
y = y(\tau) :=& q^2\prod_{m=1}^{\infty}\frac{(1-q^{2m})^2(1-q^{4m})(1-q^{5m})(1-q^{20m})^3}{(1-q^{m})^5(1-q^{10m})^2}.
\end{align}  Furthermore, we can define the associated $U$ operator sequence

\begin{align}
U^{(1-i)}\left( f \right) &:= U_5\left( \left( \frac{\mathrm{C}\Phi_2(\tau)}{\mathrm{C}\Phi_2(5\tau)} \right)^i \cdot f \right) \text{ for } i=0,1,\\
U^{(\alpha)}\left( f \right) &:= U^{(\alpha\bmod{2})}\left( f \right).
\end{align} In this case, for the function sequence $(L_{\alpha})_{\alpha\ge 1}$ associated with the Andrews--Sellers congruences, we will have

\begin{align}
L_{\alpha+1} &= U^{(\alpha)}\left( L_{\alpha} \right).
\end{align}  Let us now examine the function

\begin{align}
t&:= y+4xy = q^2 + 9q^3 + 50q^4 + 219q^5 + ...
\end{align}  Notice that $t\not\equiv 0\pmod{5}$.  One can show by simple computation that

\begin{align}
U^{(0)}\circ U^{(1)}(t)\equiv t\pmod{5}.
\end{align}  In other words, $t$ is a nontrivial eigenfunction modulo 5 with respect to the composite operator $U^{(0)}\circ U^{(1)}$.  Of course, this means that no element in the sequence

\begin{align}
\left(\left(U^{(0)}\circ U^{(1)}\right)^n (t)\right)_{n\ge 0}
\end{align} will have divisibility by 5.  Certainly, the sequence will not converge to 0 in the 5-adic sense.

The significance of this matter can be better grasped if we consider convergent sequences in the standard topology.  Suppose we want to prove that

\begin{align*}
\lim_{\alpha\rightarrow\infty} L_{\alpha} = 0
\end{align*} for a given sequence of functions

\begin{align*}
(L_{\alpha})_{\alpha\ge 0}.
\end{align*}

One possible approach is to find some reference functions, e.g. $(F_{\alpha})_{\alpha\ge 0}, (G_{\alpha})_{\alpha\ge 0}$, which we might have some better control over, such that

\begin{align*}
L_{\alpha} = F_{\alpha} + G_{\alpha}.
\end{align*}  If we wish to prove that $\lim_{\alpha\rightarrow\infty} L_{\alpha} = 0$, it is of course sufficient to show that

\begin{align*}
\lim_{\alpha\rightarrow\infty} F_{\alpha} = \lim_{\alpha\rightarrow\infty} G_{\alpha} = 0.
\end{align*}  However, it is certainly not necessary!  We could, for example, have 

\begin{align*}
\lim_{\alpha\rightarrow\infty} F_{\alpha} = 1,\ \lim_{\alpha\rightarrow\infty} G_{\alpha} = -1.
\end{align*}  If we want to prove convergence of $L_{\alpha}$ term-wise, it is clear that we need to carefully select the functions $F_{\alpha}$, $G_{\alpha}$; otherwise, such reference functions will not be very useful.  This same principle holds if we consider function sequences in the 5-adic topology.

It is precisely this complication which makes the Andrews--Sellers family so difficult to prove.  It does not yet appear possible to choose a set of natural reference functions---say, those which live at a single cusp---and to construct the necessary induction to prove the entire family, since these reference functions often turn out to be eigenfunctions themselves.

The methods developed by Paule and Radu to complete the proof are extremely powerful, and have been applied to other problems, e.g., Theorem \ref{wyfirst}.  The primary disadvantage of these methods is that, while they certainly work, a strong theoretical understanding of them is lacking.  We still do not understand what the associated algebraic structures ought to look like for any given problem.  Similarly, it is not clear how to properly select the associated reference functions for these algebraic structures.  The setup for the method is essentially an educated guess, coupled with a large degree of experimentation.

\section{Conclusion}

After more than a century, our understanding of partition congruences is at a peculiar point.  Our understanding continues to grow, but we continue to find extraordinary difficulties as new congruence families are found, all of which superficially resemble Ramanujan's original work.  A serious program to classify congruence families has only very recently begun, and a strong theoretical foundation of the most difficult of these problems has yet to be established.  There is much to be done.

We also must quickly add that we have said nothing of more complex first order congruences that were pioneered by Atkin and O'Brien (\cite{Atkin0}, \cite{Atkin3}), of which the underlying theory has been developed by Ahlgren (e.g., \cite{Ahlgren0}), Ono (e.g., \cite{Ono2}), and Radu (\cite{Radu2}).  We have also left aside the enormously rich and beautiful study of the combinatorial manifestations of partition congruences, as in the theory of ranks and cranks developed by Dyson \cite{Dyson}, Andrews, and especially Garvan (\cite{AndrewsGarv} and \cite{Garvan0}).  We have only discussed a comparatively narrow aspect of this subject.  May we continue to find more joy in the curiosity of adding whole numbers.

\section{Acknowledgments}

This research was funded in whole by the Austrian Science Fund (FWF): Einzelprojekte P 33933, ``Partition Congruences by the Localization Method".  My sincerest and humblest thanks to the Austrian Government and People for their generous support.

I wish to thank Professor Peter Paule for his guidance over the years that has shaped my understanding and perspective of this very fruitful topic.  I also want to thank Georgia Southern's Mathematics Department for their organization of a wonderful virtual conference that allowed me the pleasure of meeting old colleagues again, and provided a little relief to the rolling lockdowns of 2021.

Finally, a huge thanks to Divine Wanduku and Drew Sills for their kindness and patience with my many delays in submission.

\end{document}